# Basic Properties of Incomplete Macdonald Function with Applications

Jian-Jun SHU, and Kunal Krishnaraj SHASTRI

School of Mechanical & Aerospace Engineering, Nanyang Technological University, 50 Nanyang Avenue, Singapore 639798.

Correspondence should be addressed to Jian-Jun SHU; mjjshu@ntu.edu.sg

## Abstract

The incomplete version of the Macdonald function has various appellations in literature and earns a well-deserved reputation of being a computational challenge. This paper ties together the previously disjoint literature and presents the basic properties of the incomplete Macdonald function, such as recurrence and differential relations, series and asymptotic expansions. This paper also shows that the incomplete Macdonald function, as a simple closed-form expression, is a particular solution to a parabolic partial differential equation, which arises naturally in a wide variety of transient and diffusion-related phenomena.

## 1. Introduction

As a particular solution to the modified Bessel equation [1], the Macdonald function [2] (also known as the modified Bessel function of the second kind, the Basset function, or the modified Hankel function) has been employed in wide-ranging fields to provide analytical solutions for many physical phenomena. The Macdonald function can be expressed in the form of integrals. A new class of incomplete special function, called the incomplete Macdonald function, is defined by having the variable endpoint of integration and arises in wide-ranging contexts such as viscous flow [3], heat conduction [4], groundwater hydrology [5], electromagnetism [6], galaxy [7], nuclear reactor [8], and doubly special relativity [9].

The several distinct incomplete versions of the Macdonald function are possible since various integral representations of the Macdonald function often seem unrelated under the variable endpoint of integration. In literature, the specific definition chosen among various forms of the incomplete Macdonald function is directed by the particular application. These various forms and interrelationships of the incomplete Macdonald function are summarized in the Table 1.

Table 1: Various distinct incomplete versions of the Macdonald function

| Function | Definition |
|---|---|
| Shu function [3] | $S_\nu(z,t) = \dfrac{1}{2}\left(\dfrac{z}{2}\right)^\nu \displaystyle\int_0^t \dfrac{e^{-\tau - \frac{z^2}{4\tau}}}{\tau^{\nu+1}} d\tau$ |
| Generalized incomplete gamma function [4] | $\Gamma(\nu, t; z) = \displaystyle\int_t^\infty \tau^{\nu-1} e^{-\tau - \frac{z}{\tau}} d\tau = 2 z^{\frac{\nu}{2}} S_\nu\left(2\sqrt{z}, \dfrac{z}{t}\right)$ |
| Leaky aquifer function [5] | $L_\nu(z,t) = \displaystyle\int_1^\infty \dfrac{e^{-z\tau - \frac{t}{\tau}}}{\tau^{\nu+1}} d\tau = 2\left(\dfrac{z}{t}\right)^{\frac{\nu}{2}} S_{-\nu}\left(2\sqrt{zt}, t\right)$ |
| Incomplete modified Bessel function [6] | $\dfrac{1}{2}\displaystyle\int_t^\infty e^{-z\cosh(\tau)} \cosh(\nu\tau) d\tau = \dfrac{1}{2}(S_\nu + S_{-\nu})\left(z, \dfrac{z e^{-t}}{2}\right)$ |

The difference in the names assigned to various forms of the incomplete Macdonald function is worth noting (as listed in Table 1). They are the Shu function of viscous flow, generalized incomplete gamma function of heat conduction, leaky aquifer function of groundwater hydrology, and incomplete modified Bessel function of electromagnetism. It is pretty apparent that there is a lack of communication among different research communities. The Shu function was first employed by Shu and Chwang [3] in the expression of the hydrodynamic force acting on a rigid circular cylinder translating in a time-dependent rotating flow field. The motivation of studying the generalized incomplete gamma function [4] was the role it played in the closed-form solution to several problems in heat conduction. Groundwater hydrologists commonly refer to the leaky aquifer integral as the leaky aquifer function [5], which is useful for determining the hydraulic properties of leaky-confined aquifers. The incomplete modified Bessel function [6] was introduced to express the solution of electromagnetic problems in truncated cylindrical structures. It is worth mentioning at this end that these Shu function, generalized incomplete gamma function, leaky aquifer function, and incomplete modified Bessel function have a well-deserved reputation of being a computational challenge due to their integral representations. The increasing number of applications calls for a rigorous analysis of various forms of the incomplete Macdonald function. To avoid redundancy, the special function analyzed in this paper is mainly referred to as the Shu function, which is a representative of a large class of time-dependent problems.

In this paper, the key properties of the Shu function, such as recurrence and differential relations, series and asymptotic expansions, are derived. It is also shown that the Shu function is a particular solution to a parabolic partial differential equation (PDE), which occurs in various transient problems, for example, transient flow in porous media [10], electromagnetic waves in a cylindrical waveguide [11, 12] and diffusion-related phenomena [8].

## 2. Definition

The integral form of the Macdonald function is given [13] by

$$K_\nu(z) = \frac{1}{2}\left(\frac{z}{2}\right)^\nu \int_0^\infty \frac{e^{-\tau - \frac{z^2}{4\tau}}}{\tau^{\nu+1}} d\tau. \tag{1}$$

As an incomplete version of the Macdonald function, the definition of the Shu function [3] is adopted for investigation in this paper and obtained by restricting the upper endpoint of the integral in (1) to $t$,

$$S_\nu(z,t)=\frac{1}{2}\left(\frac{z}{2}\right)^\nu \int_0^t \frac{e^{-\tau-\frac{z^2}{4\tau}}}{\tau^{\nu+1}}d\tau, \quad \text{Re}(t)>0 \text{ and } \text{Re}(z)>0, \tag{2}$$

in which Re denotes the real part of a complex number, the symbols $\nu$, $z$, and $t$ are, respectively, termed the order, argument, and endpoint of the Shu function and are generally taken to be complex quantities. Hereafter, the symbols $n$ and $x>0$ are, respectively, annotated to denote integer order $\nu$ and real argument $z$. The behavior of the Shu function $S_n(x,t)$ of integer order $n$, real argument $x>0$, and real endpoint $t>0$ is depicted in Figures 1 and 2. In view of that the Shu function has infinite when $\text{Re}(z)=0$, we restrict attention to the complex argument $\text{Re}(z)>0$ and the complex endpoint $\text{Re}(t)>0$.

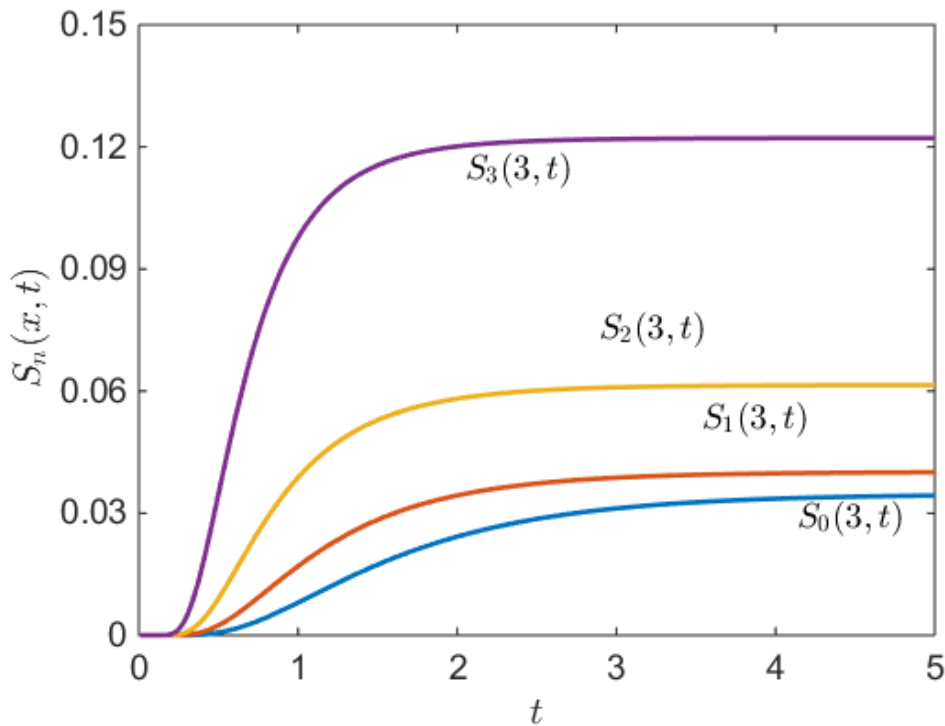

Figure 1: $S_n(x,t)$ against real endpoint $t>0$ for various values of integer order $n$ at $x=3$

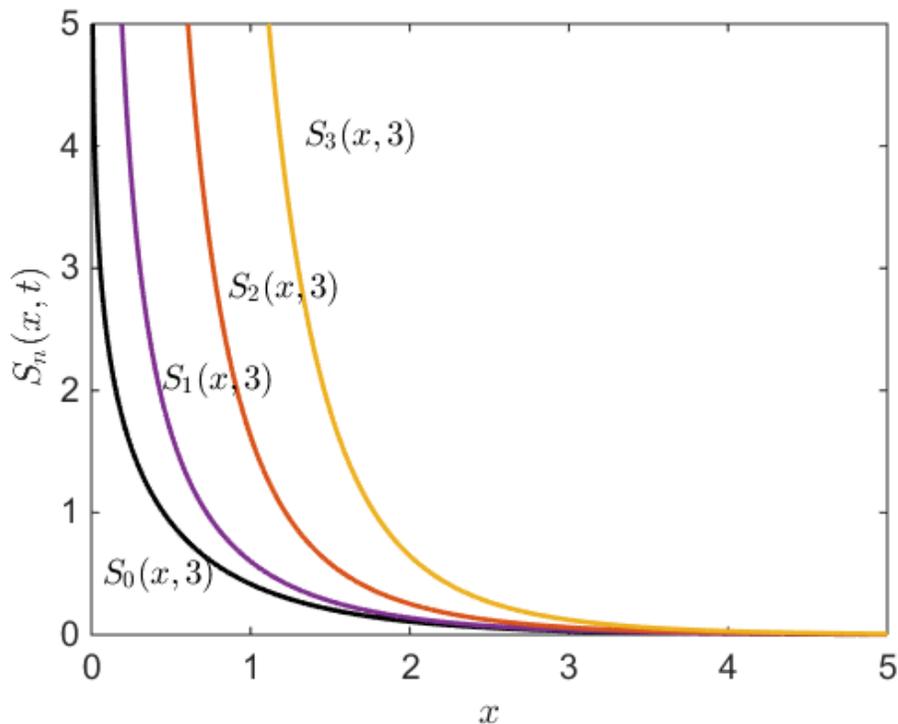

Figure 2: $S_n(x,t)$ against real argument $x>0$ for various values of integer order $n$ at $t=3$

By using (1), the Shu function can be written as

$$S_\nu(z,t) = K_\nu(z) - \frac{1}{2}\left(\frac{z}{2}\right)^\nu \int_t^\infty \frac{e^{-\tau - \frac{z^2}{4\tau}}}{\tau^{\nu+1}} d\tau. \tag{3}$$

Substituting $e^w = \frac{z}{2\tau}$ in (2), we get [14]

$$S_\nu(z,t) = \frac{1}{2} \int_{\ln\left(\frac{z}{2t}\right)}^\infty e^{-z\cosh(w)+\nu w} dw. \tag{4}$$

By using the substitution: $y = \frac{z^2}{4\tau}$ in (2), we get the alternate form

$$S_\nu(z,t) = \frac{1}{2}\left(\frac{2}{z}\right)^\nu \int_{\frac{z^2}{4t}}^\infty y^{\nu-1} e^{-y - \frac{z^2}{4y}} dy. \tag{5}$$

This can be expressed in the form of the generalized incomplete gamma function [4] (Table 1) as

$$S_\nu(z,t) = \frac{1}{2}\left(\frac{2}{z}\right)^\nu \Gamma\left(\nu, \frac{z^2}{4t}; \frac{z^2}{4}\right). \tag{6}$$

Furthermore, the Shu function can be related to the leaky aquifer function [5] (Table 1) by

$$S_\nu(z,t) = \frac{1}{2}\left(\frac{z}{2t}\right)^\nu L_{-\nu}\left(\frac{z^2}{4t}, t\right). \tag{7}$$

Using this relation, the algorithm [15] proposed for the computation of numerical value of $L_\nu(z,t)$ can be used for the computation of $S_\nu(z,t)$.

## 3. Recurrence Relations

Two recurrence relations of the Shu function are derived in this section.

### 3.1 The First Recurrence Relation

Integrating (2) by parts, we obtain

$$S_\nu(z,t) = -\frac{1}{2}\left(\frac{z}{2}\right)^\nu \left.\frac{e^{-\tau - \frac{z^2}{4\tau}}}{\nu \tau^\nu}\right|_{0^+}^t + \frac{1}{2}\left(\frac{z}{2}\right)^\nu \int_0^t \frac{e^{-\tau - \frac{z^2}{4\tau}}}{\nu \tau^\nu}\left(-1 + \frac{z^2}{4\tau^2}\right) d\tau.$$

Here, note that

$$\frac{\partial S_{\nu-1}}{\partial t} = \frac{1}{2}\left(\frac{z}{2}\right)^{\nu-1} \frac{e^{-t - \frac{z^2}{4t}}}{t^\nu}.$$

Hence, we obtain

$$S_\nu = -\frac{z}{2\nu}\frac{\partial S_{\nu-1}}{\partial t} - \frac{z}{2\nu} S_{\nu-1} + \frac{z}{2\nu} S_{\nu+1}.$$

On simplifying this expression, we obtain the first recurrence relation,

$$-\frac{2\nu}{z} S_\nu = \frac{\partial S_{\nu-1}}{\partial t} + S_{\nu-1} - S_{\nu+1}. \tag{8}$$

### 3.2 The Second Recurrence Relation

We differentiate (2) with respect to $z$,

$$\frac{\partial S_\nu}{\partial z} = \frac{\nu}{4}\left(\frac{z}{2}\right)^{\nu-1} \int_0^t \frac{e^{-\tau - \frac{z^2}{4\tau}}}{\tau^{\nu+1}} d\tau - \frac{1}{2}\left(\frac{z}{2}\right)^{\nu+1} \int_0^t \frac{e^{-\tau - \frac{z^2}{4\tau}}}{\tau^{\nu+2}} d\tau.$$

Here,

$$\int_0^t \frac{e^{-\tau - \frac{z^2}{4\tau}}}{\tau^{\nu+1}} d\tau = -\frac{e^{-t - \frac{z^2}{4t}}}{\nu t^\nu} + \int_0^t \frac{e^{-\tau - \frac{z^2}{4\tau}}}{\nu \tau^\nu}\left(-1 + \frac{z^2}{4\tau^2}\right) d\tau.$$

Hence, we obtain

$$\frac{\partial S_\nu}{\partial z} = -\frac{1}{4}\left(\frac{z}{2}\right)^{\nu-1} \frac{e^{-t - \frac{z^2}{4t}}}{t^\nu} + \frac{1}{4}\left(\frac{z}{2}\right)^{\nu-1} \int_0^t \frac{e^{-\tau - \frac{z^2}{4\tau}}}{\tau^\nu}\left(-1 + \frac{z^2}{4\tau^2}\right) d\tau - \frac{1}{2}\left(\frac{z}{2}\right)^{\nu+1} \int_0^t \frac{e^{-\tau - \frac{z^2}{4\tau}}}{\tau^{\nu+2}} d\tau.$$

On simplifying, we obtain the second recurrence relation,

$$-2\frac{\partial S_\nu}{\partial z} = \frac{\partial S_{\nu-1}}{\partial t} + S_{\nu-1} + S_{\nu+1}. \tag{9}$$

## 4. Differential Relations

Two differential relations of the Shu function are derived in this section.

### 4.1 The First Differential Relation

Adding the two recurrence relations (8) and (9), we obtain

$$-\frac{\partial S_\nu}{\partial z} - \frac{\nu}{z} S_\nu = S_{\nu-1} + \frac{\partial S_{\nu-1}}{\partial t}. \tag{10}$$

Premultiplying by $z^\nu$, we get

$$\frac{\partial}{\partial z}(z^\nu S_\nu) = -\left(1 + \frac{\partial}{\partial t}\right)(z^\nu S_{\nu-1}).$$

This differential relation may be extended as

$$\left(\frac{1}{z}\frac{\partial}{\partial z}\right)^k (z^\nu S_\nu) = (-1)^k \left(1 + \frac{\partial}{\partial t}\right)^k (z^{\nu-k} S_{\nu-k}), \quad k \in \{0, 1, 2, \cdots\}. \tag{11}$$

### 4.2 The Second Differential Relation

We subtract the second recurrence relation (9) from (8),

$$-\frac{\partial S_\nu}{\partial z} + \frac{\nu}{z} S_\nu = S_{\nu+1}. \tag{12}$$

Premultiplying by $z^{-\nu}$,

$$\frac{\partial}{\partial z}\left(\frac{S_\nu}{z^\nu}\right) = -\frac{S_{\nu+1}}{z^\nu},$$

or

$$\left(\frac{1}{z}\frac{\partial}{\partial z}\right)^k \left(\frac{S_\nu}{z^\nu}\right) = (-1)^k \left(\frac{S_{\nu+k}}{z^{\nu+k}}\right), \quad k \in \{0, 1, 2, \cdots\}. \tag{13}$$

## 5. Partial Differential Equation

Replacing $\nu$ by $\nu+1$ in (10), we obtain

$$-zS_\nu = z\frac{\partial S_\nu}{\partial t} + z\frac{\partial S_{\nu+1}}{\partial z} + (\nu+1)S_{\nu+1}. \tag{14}$$

On differentiating (12) with respect to $z$, we get

$$\nu \frac{\partial S_\nu}{\partial z} = S_{\nu+1} + z\frac{\partial S_{\nu+1}}{\partial z} + \frac{\partial S_\nu}{\partial z} + z\frac{\partial^2 S_\nu}{\partial z^2}. \tag{15}$$

From (12), (14), and (15), it is a straightforward exercise to obtain the following PDE:

$$z^2 \frac{\partial^2 S_\nu}{\partial z^2} + z\frac{\partial S_\nu}{\partial z} - (z^2 + \nu^2)S_\nu - z^2 \frac{\partial S_\nu}{\partial t} = 0.$$

The PDE can also be written as

$$\left(\Omega - \frac{\partial}{\partial t}\right)S_\nu = 0, \tag{16}$$

where $\Omega$ represents the modified Bessel operator:

$$\Omega \equiv \frac{1}{z}\frac{\partial}{\partial z}\left(z\frac{\partial}{\partial z}\right) - \left(1 + \frac{\nu^2}{z^2}\right).$$

This elegant parabolic PDE (16) arises naturally in a wide variety of transient and diffusion-related phenomena.

## 6. Series Expansions

Two expansions as $|t| \to 0$ and $|z| \to 0$, respectively, are derived in this section, in which $|\,|$ denotes the modulus of a complex number.

## 6.1 Series Expansion as $|t|\to 0$

We substitute the expansion $e^{-\frac{z^2}{4y}} = \sum_{k=0}^{\infty} \frac{(-1)^k}{k!}\left(\frac{z}{2}\right)^{2k}\frac{1}{y^k}$ as $|y|\to\infty$ in (5) to obtain

$$S_\nu(z,t) = \sum_{k=0}^{\infty} \frac{(-1)^k}{2\,k!}\left(\frac{2}{z}\right)^{\nu-2k} \int_{\frac{z^2}{4t}}^{\infty} y^{\nu-k-1} e^{-y}\,dy,$$

or

$$S_\nu(z,t) = \sum_{k=0}^{\infty} \frac{(-1)^k}{2\,k!}\left(\frac{2}{z}\right)^{\nu-2k} \Gamma\!\left(\nu-k,\frac{z^2}{4t}\right). \qquad (17)$$

The asymptotic expansion of the incomplete gamma function $\Gamma(\nu,t) = \int_t^{\infty} \tau^{\nu-1} e^{-\tau}\,d\tau$ for large $|t|$ is given [4] by

$$\Gamma(\nu,t) = t^{\nu-1} e^{-t} \sum_{m=0}^{\infty} (-1)^m (1-\nu)_m \left(\frac{1}{t}\right)^m, \qquad (18)$$

where $(\nu)_m \equiv \frac{\Gamma(m+\nu)}{\Gamma(\nu)}$ is the Pochhammer polynomial and $\Gamma(\nu) = \Gamma(\nu,0) = \int_0^{\infty} \tau^{\nu-1} e^{-\tau}\,d\tau$ is the gamma function. We substitute the expansion (17) to get

$$S_\nu(z,t) = \sum_{m=0}^{\infty}\sum_{k=0}^{\infty} \frac{(-1)^{m+k}}{2\,k!}(1-\nu+k)_m \left(\frac{z}{2}\right)^{\nu-2m-2} \frac{e^{-\frac{z^2}{4t}}}{t^{\nu-m-k-1}}.$$

The leading term approximation is given by

$$S_\nu(z,t) \sim \frac{1}{2}\left(\frac{z}{2}\right)^{\nu-2} \frac{e^{-\frac{z^2}{4t}}}{t^{\nu-1}}\left[1+O(|t|)\right] \quad \text{as } |t|\to 0. \qquad (19)$$

The leading term approximation to the actual value of the Shu function $S_n(x,t)$ of integer order $n$ and real argument $x>0$ for small endpoint $t$ is shown in Figure 3. As can be observed, the agreement is better for higher integer order $n$.

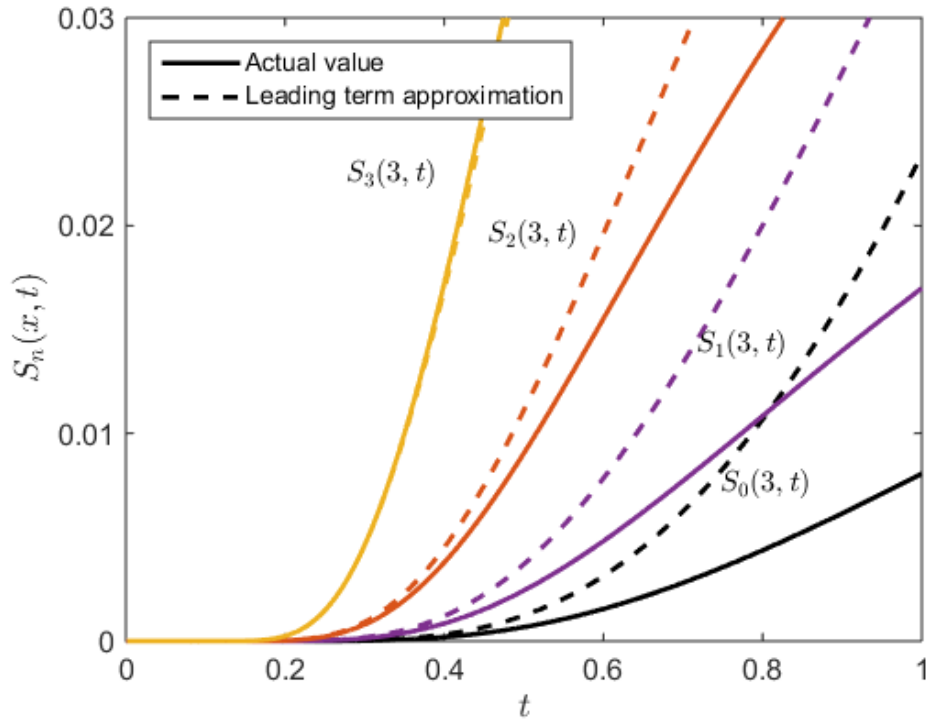

Figure 3: Leading term approximation for $t \to 0^+$

## 6.2 Series Expansion as $|z| \to 0$

We substitute the expansion $e^{-\frac{z^2}{4\tau}} = \sum_{k=0}^{\infty} \frac{(-1)^k}{k!} \frac{1}{\tau^k} \left(\frac{z}{2}\right)^{2k}$ as $|z| \to 0$ in (3) to get

$$S_\nu(z,t) = K_\nu(z) - \sum_{k=0}^{\infty} \frac{(-1)^k}{2k!} \left(\frac{z}{2}\right)^{\nu+2k} \int_t^{\infty} \frac{e^{-\tau}}{\tau^{\nu+k+1}} d\tau .$$

This expansion can be expressed in a more compact form by using the incomplete gamma function $\Gamma(\nu,t)$:

$$S_\nu(z,t) = K_\nu(z) - \sum_{k=0}^{\infty} \frac{(-1)^k \, \Gamma(-\nu-k,t)}{2k!} \left(\frac{z}{2}\right)^{\nu+2k} . \tag{20}$$

By using the approximation [4] of $K_\nu(z)$ for small $|z|$,

$$K_\nu(z) \sim \begin{cases} -\ln(z)\left[1 + O\left(\frac{1}{\ln(|z|)}\right)\right] & \nu = 0 \\ \dfrac{2^{\nu\,\mathrm{sgn}(\mathrm{Re}(\nu))-1} \Gamma(\nu\,\mathrm{sgn}(\mathrm{Re}(\nu)))}{z^{\nu\,\mathrm{sgn}(\mathrm{Re}(\nu))}} \left[1 + O\left(|z|^2 \ln(|z|)\right) + O\left(|z|^{2|\mathrm{Re}(\nu)|}\right)\right] & \nu \neq 0 \end{cases} \quad \text{as } |z| \to 0,$$

in which the signum function $\text{sgn}(y) = \begin{cases} 1 & y > 0 \\ 0 & y = 0 \\ -1 & y < 0 \end{cases} = \dfrac{2}{\pi}\int_0^\infty \dfrac{\sin(y\tau)}{\tau} d\tau$. The leading term approximation is given by

$$S_\nu(z,t) \sim \begin{cases} -\ln(z)\left[1 + O\left(\dfrac{1}{\ln(|z|)}\right)\right] & \nu = 0 \\ \dfrac{2^{\nu\,\text{sgn}(\text{Re}(\nu))-1}\,\Gamma(\nu\,\text{sgn}(\text{Re}(\nu)))}{z^{\nu\,\text{sgn}(\text{Re}(\nu))}}\left[1 + O\left(|z|^2 \ln(|z|)\right) + O\left(|z|^{2|\text{Re}(\nu)|}\right)\right] & \nu \neq 0 \end{cases}$ as $|z| \to 0$, (21)

which agrees with the known result by Shu and Chwang [3] for the case of $\nu = 0$. The leading term approximation to the actual value of the Shu function $S_n(x,t)$ of integer order $n$ and real endpoint $t > 0$ for small argument $x$ is shown in Figure 4. As can be observed, the agreement is better for lower integer order $n$.

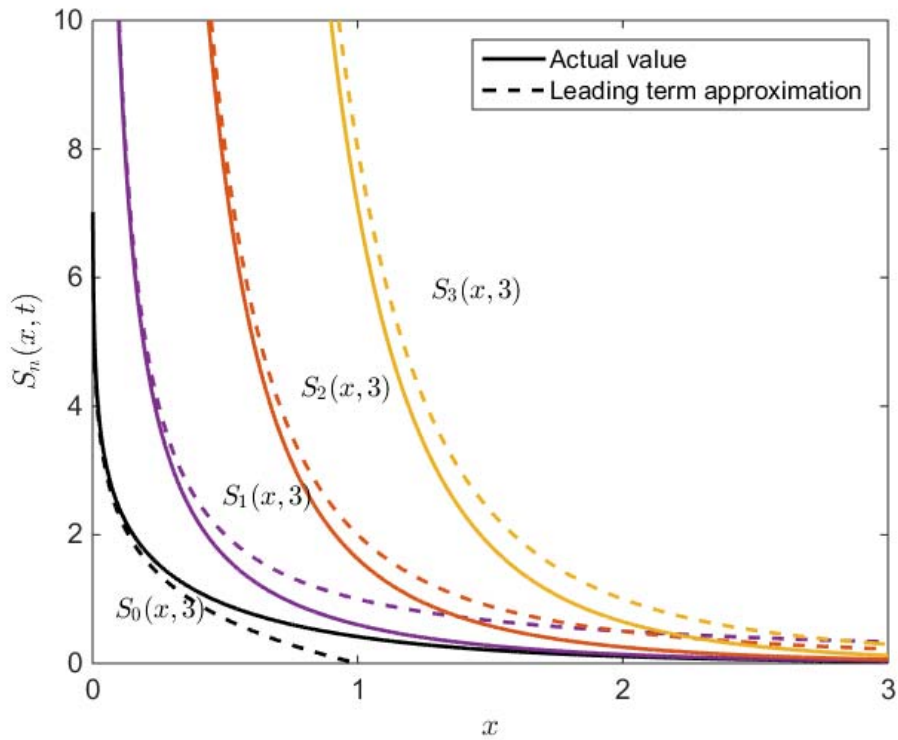

Figure 4: Leading term approximation for $x \to 0^+$

## 7. Asymptotic Expansions

Two asymptotic expansions for both large $|t|$ with fixed $z$ and large $|z|$ with fixed $t$, respectively, are derived in this section.





## 7.1 Asymptotic Expansion as $|t|\to\infty$

We use the expansions (18) and (20) to get

$$S_\nu(z,t)=K_\nu(z)-\sum_{m=0}^{\infty}\sum_{k=0}^{\infty}\frac{(-1)^{m+k}(\nu+k+1)_m}{2k!}\left(\frac{z}{2}\right)^{\nu+2k}\frac{e^{-t}}{t^{\nu+m+k+1}}.$$

The leading term approximation is given by

$$S_\nu(z,t)\sim K_\nu(z)\left[1+O\left(\frac{e^{-|\text{Re}(t)|}}{|t|^{1+\text{Re}(\nu)}}\right)\right] \quad \text{as } |t|\to\infty. \tag{22}$$

The leading term approximation to the actual value of the Shu function $S_n(x,t)$ of integer order $n$ and real argument $x>0$ for large endpoint $t$ is shown in Figure 5. As can be observed, the agreement is excellent for all integer order $n$.

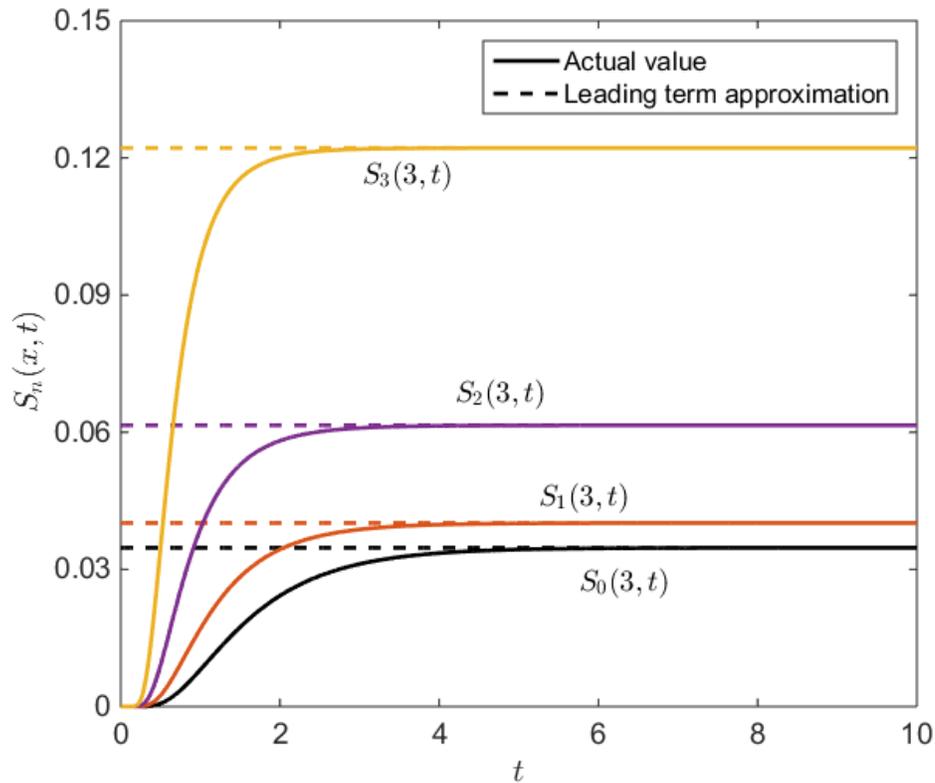

Figure 5: Leading term approximation for $t\to+\infty$





## 7.2 Asymptotic Expansion as $|z| \to \infty$

We consider an integral of the form $f(z) = \int_a^b e^{-z p(w)} q(w) dw$. If $\text{Re}(p(c)) > \text{Re}(p(a))$, $\forall c \in (a,b)$ and $q(w) = \sum_{k=0}^{\infty} q_k (w-a)^k$, $\dfrac{d\, p(w)}{d\, w} = \sum_{k=0}^{\infty} (k+1) p_k (w-a)^k$, then by the generalization of the Laplace's method [16] for complex integration,

$$f(z) = e^{-z p(a)} \sum_{k=0}^{\infty} \Gamma(k+1) \frac{c_k}{z^{k+1}}, \tag{23}$$

where the first two coefficients $c_0$ and $c_1$ are given by

$$c_0 = \frac{q_0}{p_0}, \quad c_1 = \frac{p_0 q_1 - 2 p_1 q_0}{p_0^3}.$$

We use (4) to derive the asymptotic expansion for large $|z|$ and note that $q(w) = e^{vw}$ and $p(w) = \cosh(w)$. The function $\text{Re}(p(w)) = \text{Re}(\cosh(w))$ attains minima at $\zeta = a = \ln\left(\dfrac{z}{2t}\right)$.

Hence, $q(w) = e^{vw} \sim e^{v\zeta} + v e^{v\zeta}(w-\zeta) + O\left[|w-\zeta|^2\right]$ and $\dfrac{d\, p(w)}{d\, w} = \sinh(w) \sim \sinh(\zeta) + \cosh(\zeta)(w-\zeta) + O\left[|w-\zeta|^2\right]$. Following (23), the asymptotic expansion as $|z| \to \infty$ with fixed $t$ can be easily derived:

$$S_v(z,t) \sim \frac{e^{v\zeta - z \cosh(\zeta)}}{2 z \sinh(\zeta)} \left[1 + O\left(\frac{1}{|z|}\right)\right].$$

The leading term approximation is given by

$$S_v(z,t) \sim \frac{z^v e^{-\frac{z^2}{4t} - t}}{(2t)^{v-1}(z^2 - 4t^2)} \left[1 + O\left(\frac{1}{|z|}\right)\right] \quad \text{as } |z| \to \infty. \tag{24}$$

It is worth mentioning at this end that the incomplete modified Bessel function,

$$\frac{1}{2} \int_t^{\infty} e^{-z \cosh(\tau)} \cosh(v\tau) d\tau = \frac{1}{2}(S_v + S_{-v})\left(z, \frac{z e^{-t}}{2}\right) \sim \frac{\cosh(vt) e^{-z \cosh(t)}}{2 z \sinh(t)} \left[1 + O\left(\frac{1}{|z|}\right)\right] \quad \text{as } |z| \to \infty,$$

which agrees with the result of Cicchetti and Faraone [6], but their expression is extremely complicated. The leading term approximation to the actual value of the Shu function $S_n(x,t)$ of integer order $n$ and real endpoint $t > 0$ for large argument $x$ is shown in Figure 6. As can be observed, the agreement is better for lower integer order $n$.

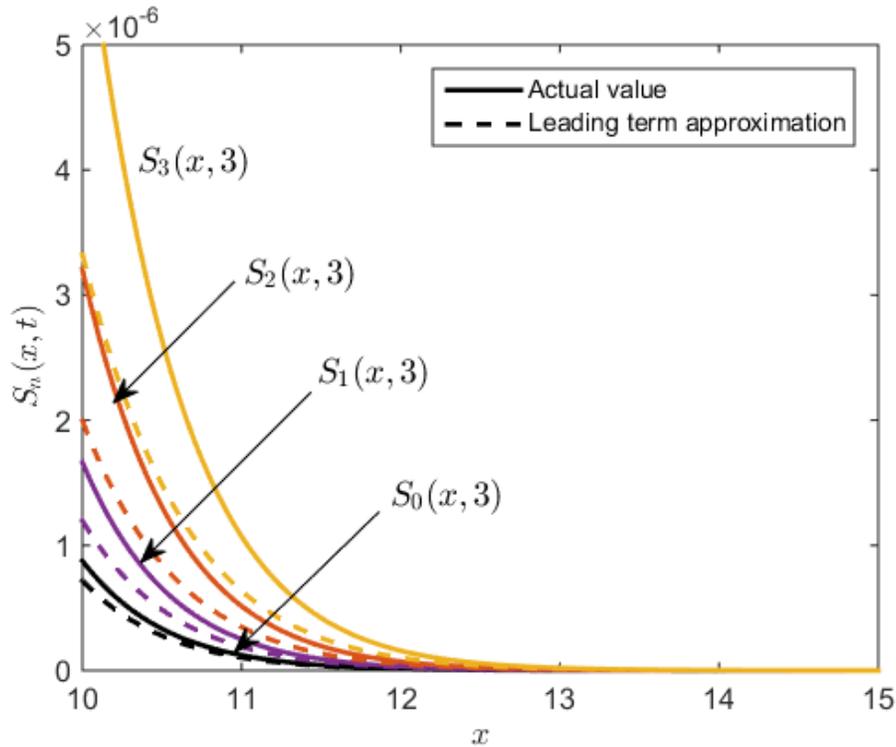

Figure 6: Leading term approximation for $x \to +\infty$

## 8. Conclusions

A number of key properties of the incomplete Macdonald function have been derived, comprising recurrence and differential relations, series and asymptotic expansions. As can be observed graphically, the agreement is found between the actual value and the corresponding leading term approximation for all limiting cases. It is also shown that the incomplete Macdonald function is a particular solution to the parabolic PDE associated with a wide variety of transient natural phenomena. The Shu function and other distinct incomplete versions of the Macdonald function can be used to find simple closed-form expressions for various natural phenomena.

## Data Availability

No data were used to support this study.

## Conflicts of Interest

The authors declare that there is no conflict of interest regarding the publication of this paper.

## Acknowledgments

This work was supported by Singapore Ministry of Education Academic Research Fund Tier 1 (04MNP002133C160).

Source: Journal of Function Spaces, Vol. 2020, pp. 6548298, 2020;
DOI: 10.1155/2020/6548298## References

[1]  F.W. Bessel, "Analytische auflösung der keplerschen aufgabe," *Abhandlungen der Mathematischen Klasse der Königlich Preußischen Akademie der Wissenschaften zu Berlin*, vol. 25, pp. 49–55, 1819. (in German)

[2]  H.M. Macdonald, "Zeroes of the Bessel functions," *Proceedings London Mathematical Society*, vol. s1-30, no. 1, pp. 165–179, 1898.

[3]  J.-J. Shu and A.T. Chwang, "Generalized fundamental solutions for unsteady viscous flows," *Physical Review E*, vol. 63, no. 5, pp. 051201, 2001.

[4]  M.A. Chaudhry and S.M. Zubair, "Generalized incomplete gamma functions with applications," *Journal of Computational and Applied Mathematics*, vol. 55, no. 1, pp. 99–123, 1994.

[5]  F.E. Harris, "On Kryachko's formula for the leaky aquifer function," *International Journal of Quantum Chemistry*, vol. 81, no. 5, pp. 332–334, 2001.

[6]  R. Cicchetti and A. Faraone, "Incomplete Hankel and modified Bessel functions: A class of special functions for electromagnetics," *IEEE Transactions on Antennas and Propagation*, vol. 52, no. 12, pp. 3373–3389, 2004.

[7]  D.N. Spergel, "Analytical galaxy profiles for photometric and lensing analysis," *Astrophysical Journal Supplement Series*, vol. 191, no. 1, pp. 58–65, 2010.

[8]  M.M. Agrest and M.S. Maksimov, *Theory of Incomplete Cylindrical Functions and their Applications*, Springer, 2011.

[9]  N. Chandra and S. Chatterjee, "Thermodynamics of ideal gas in doubly special relativity," *Physical Review D*, vol. 85, no. 4, pp. 045012, 2012.

[10]  M.S. Hantush, "Analysis of data from pumping tests in leaky aquifers," *Transactions of the American Geophysical Union*, vol. 37, no. 6, pp. 702–714, 1956.

[11]  R. Cicchetti, A. Faraone, G. Orlandi and D. Caratelli, "Real-argument incomplete Hankel functions: Accurate and computationally efficient integral representations and their asymptotic approximants," *IEEE Transactions on Antennas and Propagation*, vol. 63, no. 6, pp. 2751–2756, 2015.

[12]  R. Cicchetti, V. Cicchetti, A. Faraone and O. Testa, "Analysis of thin truncated cylinder scatterers using incomplete Hankel functions and surface impedance boundary conditions," *IEEE Access*, vol. 8, pp. 72997–73004, 2020.

[13]  G.N. Watson, *A Treatise on the Theory of Bessel Functions*, Nabu Press, 2014.

[14]  J.-J. Shu and J.S. Lee, "Fundamental solutions for micropolar fluids," *Journal of Engineering Mathematics*, vol. 61, no. 1, pp. 69–79, 2008.

[15]  R.M. Slevinsky and H. Safouhi, "A recursive algorithm for the G transformation and accurate computation of incomplete Bessel functions," *Applied Numerical Mathematics*, vol. 60, no. 12, pp. 1411–1417, 2010.

[16]  P.S. Laplace, "Mémoire sur la probabilité des causes sur les évènements," *Mémoires de Mathématique et de Physique, Tome Sixième*, 1774. (in French)14